\documentclass[11pt,reqno]{amsart}
\usepackage{etoolbox}
\makeatletter
\patchcmd{\@settitle}{\uppercasenonmath\@title}{}{}
  {\PackageError{titlecase}{Unable to patch the main title}{}}
\patchcmd{\maketitle}{\uppercasenonmath\shorttitle}{}{}
  {\PackageError{titlecase}{Unable to patch the running title}{}}
\patchcmd{\ps@headings}{\rlap{\thepage}}{}{}
  {\PackageError{pagestyle}{Unable to patch the even head}{}}
\patchcmd{\ps@headings}{\llap{\thepage}}{}{}
  {\PackageError{pagestyle}{Unable to patch the odd head}{}}
\apptocmd{\ps@headings}{%
  \def\@oddfoot{\normalfont\scriptsize\hfil\thepage\hfil}%
  \let\@evenfoot\@oddfoot}{}
  {\PackageError{pagestyle}{Unable to patch the foot}{}}
\makeatother

\usepackage[a4paper,textwidth=15.5cm,textheight=23cm,centering]{geometry}
\usepackage{amsmath}
\usepackage{amssymb}
\usepackage[numbers,sort&compress]{natbib}
\usepackage{hyperref}
\hypersetup{colorlinks=true,linkcolor=blue,filecolor=blue,urlcolor=blue,citecolor=cyan}

\newtheorem{thm}{Theorem}[section]
\newtheorem{pro}[thm]{Proposition}
\newtheorem{lem}[thm]{Lemma}
\newtheorem{cor}[thm]{Corollary}

\newtheoremstyle{italicremark}
{}{}{\itshape}{}{\itshape}{.}{0.5em}{}
\theoremstyle{italicremark}
\newtheorem*{rem}{Remark}
\newtheorem{defi}[thm]{Definition}

\begin{document}

\title[Cayley Regression for Finite Groups with a Cyclic Subgroup of Index Two]{Cayley Regression for Finite Groups with a Cyclic Subgroup of Index Two}

\author[Y.B. Jiang] {Yingbi Jiang$^{a}$}
\author[P.f. Guo] {Pengfei Guo$^{a}$}
\author[J.H. Xie] {Jin-Hua Xie$^{a,*}$}

\begingroup
\renewcommand{\thefootnote}{}
\begin{NoHyper}
\footnotetext{
	a: School of Mathematics and Statistics, Hainan Normal University,
	Haikou, Hainan, P.R. China.
	Email: yingbijiang2727@163.com (Yingbi Jiang);
	guopf999@163.com (Pengfei Guo);
	jinhuaxie@hainnu.edu.cn (Jin-Hua Xie).
	$*$: Corresponding author.
}
\end{NoHyper}
\endgroup

	\begin{abstract}
A finite group $G$ is a Cayley regression if every generalized Cayley graph of $G$ is isomorphic to a Cayley graph of $G$, and a $k$\nobreakdash-Cayley regression if this holds for all such graphs of valency at most $k$.
We show that the dihedral group $D_{2n}$, $n\geq2$, is a $2$\nobreakdash-Cayley regression if and only if it is a $3$\nobreakdash-Cayley regression, which holds exactly when $n$ is a prime power, and that for $k\geq4$ it is a $k$\nobreakdash-Cayley regression if and only if $n$ is an odd prime power or $n\in\{2,4\}$.
Groups with a cyclic subgroup of index two that are neither cyclic nor dihedral are not even $2$\nobreakdash-Cayley regressions; consequently, a group with a cyclic subgroup of index two is a Cayley regression if and only if it is isomorphic to $C_2$, $C_4$, $C_6$, $C_{10}$, or $D_{2n}$ with $n$ as above.
Finally, every split metacyclic group $C_m\rtimes C_q$ with $m>2$ and odd $q\geq3$ is not a $2$\nobreakdash-Cayley regression.
\end{abstract}

\maketitle
\section{Introduction}

All groups and graphs considered in this paper are finite, and all graphs are simple and undirected.
For standard graph-theoretic terminology, we refer to
Godsil and Royle~\cite{GodsilRoyle2001}.

Let $G$ be a group and let $S\subseteq G\setminus\{1\}$ with $S=S^{-1}$, where $X^{-1}=\{x^{-1}\mid x\in X\}$ for any subset $X$ of $G$.
The \emph{Cayley graph} $\operatorname{Cay}(G,S)$ is the graph with vertex set $G$, in which $g$ and $h$ are adjacent if and only if $g^{-1}h\in S$.
Generalized Cayley graphs, introduced by Maru\v{s}i\v{c} et al.~\cite{MarusicScapellatoSalvi1992} in their study of canonical double covers of graphs, extend this construction by twisting the adjacency rule with an automorphism of $G$.

\begin{defi}\label{def:generalized-Cayley}
Let $G$ be a group, $S\subseteq G$, and $\alpha\in\operatorname{Aut}(G)$ such that
\begin{enumerate}
	\item\label{gc:i} $\alpha^2=\mathrm{id}$;
	\item\label{gc:ii} $\alpha(g^{-1})g\notin S$ for every $g\in G$;
	\item\label{gc:iii} for all $g,h\in G$, if
		$\alpha(g^{-1})h\in S$, then
		$\alpha(h^{-1})g\in S$.
\end{enumerate}
The \emph{generalized Cayley graph} $GC(G,S,\alpha)$ is the graph with vertex set $G$, in which $g$ and $h$ are adjacent if and only if $\alpha(g^{-1})h\in S$.
The set $S$ is called a \emph{generalized Cayley subset} of $G$ induced by $\alpha$.
\end{defi}

Put $\omega_\alpha(G)=\{\alpha(g^{-1})g\mid g\in G\}$.
Conditions~\ref{gc:ii} and~\ref{gc:iii} of Definition~\ref{def:generalized-Cayley} are equivalent to
$S\cap\omega_\alpha(G)=\varnothing$ and
$\alpha(S)=S^{-1}$, respectively; they ensure that
$GC(G,S,\alpha)$ has no loops and that its adjacency relation is symmetric.
The neighborhood of a vertex $g$ is $\alpha(g)S$, so $GC(G,S,\alpha)$ is $|S|$-valent.
When $\alpha=\mathrm{id}$, we have $\omega_\alpha(G)=\{1\}$, so conditions~\ref{gc:ii} and~\ref{gc:iii} of Definition~\ref{def:generalized-Cayley} become $1\notin S$ and $S=S^{-1}$, and $GC(G,S,\mathrm{id})$ is precisely the Cayley graph $\operatorname{Cay}(G,S)$.

Every Cayley graph is vertex-transitive, but a generalized Cayley graph need not be a Cayley graph, even when it is vertex-transitive.
Indeed, Hujdurovi\'{c}, Kutnar and Maru\v{s}i\v{c}~\cite{HujdurovicKutnarMarusic2015} constructed infinite families of vertex-transitive generalized Cayley graphs that are not Cayley graphs.
Hujdurovi\'{c} et al.~\cite{HujdurovicEtAl2017} gave generalized Cayley graphs that are not vertex-transitive, studied automorphisms and structural properties of generalized Cayley graphs, and proved that every generalized Cayley graph of order twice a prime is a Cayley graph.
The connection with canonical double covers was further studied by Hujdurovi\'{c}~\cite{Hujdurovic2019}.

For an abelian group written additively, taking $\alpha$ to be the inversion $g\mapsto-g$ makes $g$ and $h$ adjacent if and only if $g+h\in S$; thus the resulting generalized Cayley graphs are loopless Cayley sum graphs, also known as addition Cayley graphs.
Grynkiewicz et al.~\cite{GrynkiewiczLevSerra2009} determined the vertex connectivity of addition Cayley graphs on finite abelian groups.
For generalized Cayley graphs, Zhu et al.~\cite{ZhuYangChen2023} investigated their spectra on finite abelian groups, and
Liao et al.~\cite{LiaoLiuZhang2025} studied bipartiteness, connectedness and cubic integral examples.

The isomorphism problem offers another perspective.
For Cayley graphs, we refer to Babai~\cite{Babai1977}, Muzychuk~\cite{Muzychuk2004} and the survey of Li~\cite{Li2002}.
For generalized Cayley graphs, isomorphism properties were studied by Yang et al.\ for alternating groups~\cite{YangLiuChenFeng2017} and in a general setting~\cite{YangLiuFeng2018},
by Liao and Liu~\cite{LiaoLiu2023} for cyclic, dihedral and other groups, and by Zhu et al.~\cite{ZhuLiuYang2023} for finite non-abelian simple groups.
To study when a generalized Cayley graph of a group is isomorphic to a Cayley graph of the same group, Yang et al.~\cite{YangLiuFeng2018} introduced the following notion.

\begin{defi}\label{def:Cayley-regression}
	Let $G$ be a group.
	\begin{enumerate}
		\item\label{cr:i} $G$ is a \emph{Cayley regression}
		if every generalized Cayley graph of $G$ is isomorphic to a Cayley graph of $G$.
		\item\label{cr:ii} For a positive integer $k$, $G$ is a \emph{$k$\nobreakdash-Cayley regression} if every generalized Cayley graph of $G$ of valency at most $k$ is
		isomorphic to a Cayley graph of $G$.
	\end{enumerate}
\end{defi}

\begin{rem}
	Since $GC(G,S,\mathrm{id})=\operatorname{Cay}(G,S)$,
	only nonidentity automorphisms matter in Definition~\ref{def:Cayley-regression}.
	Hence, throughout this paper, we assume that
	$\alpha\neq\mathrm{id}$ unless otherwise stated.
\end{rem}

Yang et al.~\cite[Theorem~7.5]{YangLiuFeng2018}
proved that a nontrivial cyclic group of odd order is a Cayley regression if and only if its order is a prime power,
and asked for a classification of Cayley regressions among cyclic and dihedral groups
\cite[Question~7.11]{YangLiuFeng2018}.
Liao et al.~\cite[Theorem~2.11]{LiaoLiuTang} completed the cyclic case: a nontrivial cyclic group is a Cayley regression if and only if its order is an odd prime power or one of $2,4,6,10$.
They also proved that $D_{2p^m}$ is a $3$\nobreakdash-Cayley regression for every prime $p$ with $p^m\geq3$.

In this paper, we settle the dihedral case of \cite[Question~7.11]{YangLiuFeng2018} and extend the classification to all groups with a cyclic subgroup of index two.
For dihedral groups, we determine the $k$\nobreakdash-Cayley regressions for every $k\geq2$, and hence the Cayley regressions.

\begin{thm}\label{dihedral-2-cayley-regression}
	Let $n\geq2$. Then the following are equivalent:
	\begin{enumerate}
		\item\label{d23:i} the dihedral group $D_{2n}$ is a $2$\nobreakdash-Cayley regression;
		\item\label{d23:ii} the dihedral group $D_{2n}$ is a $3$\nobreakdash-Cayley regression;
		\item\label{d23:iii} $n$ is a prime power.
	\end{enumerate}
\end{thm}

\begin{thm}\label{thm:dihedral-k-cayley-regression}
	Let $n\geq2$ and $k\geq4$. Then the dihedral group $D_{2n}$ is a
	$k$\nobreakdash-Cayley regression if and only if $n$ is an odd prime
	power or $n\in\{2,4\}$.
\end{thm}

Comparing Theorems~\ref{dihedral-2-cayley-regression} and~\ref{thm:dihedral-k-cayley-regression},
we see that for $n=2^m$ with $m\geq3$, the group $D_{2n}$ is a $3$\nobreakdash-Cayley regression
but not a $4$\nobreakdash-Cayley regression.
Since every generalized Cayley graph of $D_{2n}$ has valency at most $2n-1$,
taking $k\geq\max\{4,2n-1\}$ in Theorem~\ref{thm:dihedral-k-cayley-regression} yields the classification of dihedral Cayley regressions.

\begin{cor}\label{classification-Cayley-regression}
	Let $n\geq2$. Then the dihedral group $D_{2n}$ is a
	Cayley regression if and only if $n$ is an odd prime
	power or $n\in\{2,4\}$.
\end{cor}

Groups with a cyclic subgroup of index two that are neither cyclic nor dihedral, such as the dicyclic groups $T_{4n}$ and the groups $SD_{8n}$ (see Corollary~\ref{cor:T4n}), are not even $2$\nobreakdash-Cayley regressions.
Combined with Corollary~\ref{classification-Cayley-regression} and the cyclic case, this gives the following classification.

\begin{thm}\label{thm:index-two-classification}
	Let $G$ be a group with a cyclic subgroup of index two.
	\begin{enumerate}
		\item\label{idx:i} If $G$ is neither cyclic nor dihedral, then $G$ is not a $2$\nobreakdash-Cayley regression.
		\item\label{idx:ii} $G$ is a Cayley regression if and only if $G$ is isomorphic
		to $C_2$, $C_4$, $C_6$, $C_{10}$, or $D_{2n}$, where $n$ is
		an odd prime power or $n\in\{2,4\}$.
	\end{enumerate}
\end{thm}

Finally, replacing index two by an odd index, we show that split metacyclic groups behave like the groups in Theorem~\ref{thm:index-two-classification}\ref{idx:i}.

\begin{thm}\label{thm:split-metacyclic-odd}
	Let $G=C_m\rtimes C_q$ with an arbitrary action of $C_q$ on $C_m$,
	where $m>2$ and $q\geq3$ is odd.
	Then $G$ is not a $2$\nobreakdash-Cayley regression.
\end{thm}

After collecting preliminaries in Section~\ref{Preliminaries},
we prove Theorems~\ref{dihedral-2-cayley-regression}
and~\ref{thm:dihedral-k-cayley-regression}, together with Corollary~\ref{classification-Cayley-regression}, in Section~\ref{sec:dihedral},
and Theorems~\ref{thm:index-two-classification} and~\ref{thm:split-metacyclic-odd}
in Section~\ref{sec:index-two}.

\section{Preliminaries}\label{Preliminaries}

We first fix notation.
For disjoint sets $X$ and $Y$, we write $X\mathbin{\dot\cup}Y$ for their union.
For a positive integer $n$, we write $\mathbb Z_n$ for the ring of integers modulo $n$
and $\mathbb Z_n^{*}$ for its group of units.
For a group $G$, we write $\operatorname{Aut}(G)$ for its automorphism group and
$\mathrm{id}$ for the identity automorphism.
The symbol $C_n$ denotes a cyclic group of order $n$ when used in a group expression,
and a cycle of length $n$ when used in a graph expression.
Throughout, $D_{2n}=\langle a,b\mid a^n=b^2=1,\ bab=a^{-1}\rangle$, where $n\geq2$,
denotes the dihedral group of order $2n$.

For a graph $\Gamma$ and $X\subseteq V(\Gamma)$, let $\Gamma[X]$ denote the subgraph
induced by $X$, and let $\overline{\Gamma}$ denote the complement of $\Gamma$.
The symbols $K_n$ and $K_{r,s}$ denote the complete graph and the complete bipartite graph,
and $\Gamma_1\square\Gamma_2$ denotes the Cartesian product of graphs.
For a positive integer $k$, $k\Gamma$ denotes the disjoint union of $k$ copies of $\Gamma$.
For a group $G$ and $T\subseteq G\setminus\{1\}$ with $T=T^{-1}$, the Cayley graph $\operatorname{Cay}(G,T)$
is connected if and only if $\langle T\rangle=G$, and
$\overline{\operatorname{Cay}(G,T)}=\operatorname{Cay}(G,G\setminus(\{1\}\cup T))$;
in particular, the complement of a Cayley graph of $G$ is again a Cayley graph of $G$.

The following isomorphism property allows us to replace an automorphism by any of its conjugates.

\begin{pro}[\cite{HujdurovicEtAl2017}, Proposition~2.1]\label{prop:conjugate-GC}
Let $S$ be a generalized Cayley subset of $G$ induced by
$\alpha$. 
For every $\beta\in\operatorname{Aut}(G)$, $\beta(S)$ is a generalized Cayley subset of $G$
induced by $\beta\alpha\beta^{-1}$, and the map $g\mapsto\beta(g)$ is an isomorphism from $GC(G,S,\alpha)$
to $GC(G,\beta(S),\beta\alpha\beta^{-1})$.
\end{pro}

Consequently, when studying generalized Cayley graphs up to isomorphism, it suffices to
consider one representative from each conjugacy class of involutory automorphisms in $\operatorname{Aut}(G)$.

We next record a restriction on groups admitting connected $2$-valent Cayley graphs.
This fact is well known; we include a proof for completeness.

\begin{lem}\label{lem:connected-2-valent-Cayley}
	If a group $G$ admits a connected $2$-valent Cayley graph,
	then $G$ is either cyclic or dihedral.
\end{lem}

\begin{proof}
Let $\operatorname{Cay}(G,T)$ be a connected $2$-valent Cayley graph.
Then $|T|=2$, $T=T^{-1}$ and $\langle T\rangle=G$.
If $T$ contains an element $x$ of order greater than two, then
$T=\{x,x^{-1}\}$ and $G=\langle x\rangle$ is cyclic.

Otherwise, $T=\{x,y\}$ consists of two distinct involutions.
Put $c=xy$ and let $n$ be the order of $c$; since $x\neq y$, we have $n\geq2$.
Then $y=xc$, so $G=\langle c,x\rangle$, and $xcx=yx=c^{-1}$.
Thus $x$ normalizes $\langle c\rangle$, and hence
$G=\langle c\rangle\cup\langle c\rangle x$.
If $x\in\langle c\rangle$, then also $y=xc\in\langle c\rangle$, so the cyclic group
$\langle c\rangle$ would contain two distinct involutions, which is impossible.
Hence $x\notin\langle c\rangle$ and $|G|=2n$.
As $c$ and $x$ satisfy the defining relations of $D_{2n}$ and $|G|=2n$, we get $G\cong D_{2n}$.
\end{proof}

Finally, we record two facts on dihedral groups: a result of Liao et al.\ on
generalized Cayley graphs of valency at most three, and the case of $D_4$.

\begin{lem}[\cite{LiaoLiuTang}, Theorem~3.14]
	\label{lem:dihedral-prime-power-3}
	Let $p$ be a prime and let $m\geq1$ with $p^m\geq3$.
	Then $D_{2p^m}$ is a $3$\nobreakdash-Cayley regression.
\end{lem}

\begin{lem}\label{lem:D4-Cayley}
	The group $D_4\cong C_2\times C_2$ is a Cayley regression.
\end{lem}

\begin{proof}
Let $G=D_4=\{1,a,b,ab\}$.
Every generalized Cayley graph $GC(G,S,\alpha)$ of $G$ is an $|S|$-valent graph
on four vertices, where $0\leq|S|\leq3$.
A $2$-valent graph is a disjoint union of cycles of length at least three,
so on four vertices it must be $C_4$.
Hence, up to isomorphism, the $k$-valent graphs on four vertices for
$k=0,1,2,3$ are $4K_1$, $2K_2$, $C_4$ and $K_4$, respectively.
Since these are the Cayley graphs $\operatorname{Cay}(G,\varnothing)$,
$\operatorname{Cay}(G,\{a\})$, $\operatorname{Cay}(G,\{a,b\})$ and
$\operatorname{Cay}(G,G\setminus\{1\})$, the group $D_4$ is a Cayley regression.
\end{proof}

\section{Cayley regressions for dihedral groups}\label{sec:dihedral}

In this section, we classify the dihedral groups that are $k$\nobreakdash-Cayley regressions for $k\geq2$, and those that are Cayley regressions.
Throughout this section, $a$ and $b$ are the generators of $D_{2n}$ from Section~\ref{Preliminaries}, and $A=\langle a\rangle$.
Every element of $D_{2n}$ can be written uniquely as $a^k$ or $a^kb$ with $k\in\mathbb Z_n$,
and every element of $Ab$ is an involution.
We begin by excluding the cases in which $n$ is not a prime power.

\begin{lem}\label{lem:not-prime-power}
	Let $G=D_{2n}$, where $n\geq2$ is not a prime power.
	Then $G$ is not a $2$\nobreakdash-Cayley regression.
\end{lem}

\begin{proof}
In each of the two cases below, we choose an involutory automorphism
$\alpha$ of $G$ with $\alpha(A)=A$ and a generalized Cayley subset
$S\subseteq A$ of size two.
Then $\alpha(Ab)=Ab$, so $\Gamma=GC(G,S,\alpha)$ has no edges between
$A$ and $Ab$.
We show that $\Gamma$ contains two connected components of different orders.
Since all connected components of a vertex-transitive graph are isomorphic,
$\Gamma$ is then not vertex-transitive, and hence not a Cayley graph.

\emph{Case 1: $n$ is odd.}
Write $n=Pq$, where $P>1$ is a power of an odd prime, $q>1$, and $\gcd(P,q)=1$;
since $n$ is odd, so is $q$.
By the Chinese Remainder Theorem, there exists $u\in\mathbb Z_n$ with
$u\equiv1\pmod P$ and $u\equiv-1\pmod q$.
Then $u^2\equiv1\pmod n$, so $\alpha(a)=a^u$ and $\alpha(b)=b$ define an
automorphism $\alpha$ of $G$ with $\alpha^2=\mathrm{id}$.
For $g=a^k$ and $g=a^kb$, we have $\alpha(g^{-1})g=a^{(1-u)k}$ and
$\alpha(g^{-1})g=a^{(u-1)k}$, respectively.
Since $u\equiv1\pmod P$, the integer $1-u$ is divisible by $P$;
since $1-u\equiv2\pmod q$ and $q$ is odd, it is coprime to $q$.
Hence $\gcd(1-u,n)=P$, and therefore $\omega_\alpha(G)=\langle a^P\rangle$.

Let $S=\{a,a^{-u}\}$.
As $-u\equiv-1\pmod P$ and $P$ is odd, neither $1$ nor $-u$ is divisible by $P$,
and $1\not\equiv-u\pmod n$.
Thus $|S|=2$ and $S\cap\omega_\alpha(G)=\varnothing$.
Moreover, $\alpha(a)=a^u=(a^{-u})^{-1}$ and $\alpha(a^{-u})=a^{-u^2}=a^{-1}$,
so $\alpha(S)=S^{-1}$.
The neighbors of $a^k$ are $a^{uk+1}$ and $a^{u(k-1)}$.
Hence, if $f$ denotes the permutation $k\mapsto uk+1$ of $\mathbb Z_n$,
whose inverse is $k\mapsto u(k-1)$, then the connected components of
$\Gamma[A]$ correspond to the orbits of $\langle f\rangle$ on $\mathbb Z_n$.
Since $u\equiv1\pmod P$ and $u\equiv-1\pmod q$, we have
\[
f(k)\equiv k+1\pmod P\quad\text{and}\quad f(k)\equiv1-k\pmod q
\]
for every $k\in\mathbb Z_n$.
Iterating, we obtain $f^\ell(k)\equiv k+\ell\pmod P$, while
$f^\ell(k)\equiv k\pmod q$ if $\ell$ is even and $f^\ell(k)\equiv1-k\pmod q$ if $\ell$ is odd.
Indeed, each application of $f$ adds $1$ modulo $P$, while applying $k\mapsto1-k$ twice
gives $1-(1-k)=k$, so $f^2$ is the identity modulo $q$.
By the Chinese Remainder Theorem, $f^\ell(k)=k$ in $\mathbb Z_n$ if and only if
both congruences $f^\ell(k)\equiv k\pmod P$ and $f^\ell(k)\equiv k\pmod q$ hold.
Consequently, $f^\ell(k)=k$ if and only if $\ell$ is divisible by $P$ and either $\ell$ is even
or $1-k\equiv k$, that is, $2k\equiv1\pmod q$.
Since $q$ is odd, we may choose $k_0,k_1\in\mathbb Z_n$ with $2k_0\equiv1\pmod q$ and
$2k_1\not\equiv1\pmod q$.
For $k_0$, the least $\ell\geq1$ with $f^\ell(k_0)=k_0$ is $\ell=P$, so the orbit of $k_0$
has length $P$.
For $k_1$, the integer $\ell$ must be even and divisible by $P$; since $P$ is odd, this means
that $\ell$ is divisible by $2P$, so the orbit of $k_1$ has length $2P$.
Hence the connected components of $\Gamma$ containing $a^{k_0}$ and $a^{k_1}$ have orders
$P$ and $2P$, respectively, which are different.

\emph{Case 2: $n$ is even.}
Write $n=2^et$, where $e\geq1$ and $t>1$ is odd.
Since $n$ is even, $\alpha(a)=a$ and $\alpha(b)=a^{n/2}b$ define an
automorphism $\alpha$ of $G$ with $\alpha^2=\mathrm{id}$, and a direct computation
gives $\omega_\alpha(G)=\{1,a^{n/2}\}$.
Let $i=2^{e-1}$ and $S=\{a^i,a^{-i}\}$.
Since $t>1$, we have $i\not\equiv0,n/2$ and $2i\not\equiv0\pmod n$.
Hence $|S|=2$ and $S\cap\omega_\alpha(G)=\varnothing$; clearly
$\alpha(S)=S=S^{-1}$.
The neighbors of $a^k$ are $a^{k\pm i}$, so the connected components of
$\Gamma[A]$ are the cosets of $\langle a^i\rangle$, of order
$n/\gcd(n,i)=2t$.
The neighbors of $a^kb$ are
$\alpha(a^kb)a^{\pm i}=a^{k+n/2\mp i}b$.
Put $d=n/2-i=2^{e-1}(t-1)$; since $n/2+i\equiv-d\pmod n$, the connected
components of $\Gamma[Ab]$ are the sets $\{a^{k+\ell d}b\mid\ell\in\mathbb Z\}$.
As $t$ is odd, $\gcd(n,d)=2^{e-1}\gcd(2t,t-1)=2^e$, so these components
have order $n/2^e=t$, which differs from $2t$.

Hence, in both cases, $G$ is not a $2$\nobreakdash-Cayley regression.
\end{proof}

Combining Lemma~\ref{lem:not-prime-power} with the known
prime-power case (Lemma~\ref{lem:dihedral-prime-power-3}),
we now prove Theorem~\ref{dihedral-2-cayley-regression}.

\begin{proof}[Proof of Theorem~\ref{dihedral-2-cayley-regression}]
The implication \ref{d23:ii}$\Rightarrow$\ref{d23:i} follows from the definitions, and
\ref{d23:i}$\Rightarrow$\ref{d23:iii} is Lemma~\ref{lem:not-prime-power}.
For \ref{d23:iii}$\Rightarrow$\ref{d23:ii}, let $n$ be a prime power.
If $n=2$, then $D_4$ is a Cayley regression by
Lemma~\ref{lem:D4-Cayley}, and in particular a $3$\nobreakdash-Cayley regression;
if $n\geq3$, then $D_{2n}$ is a $3$\nobreakdash-Cayley regression by
Lemma~\ref{lem:dihedral-prime-power-3}.
\end{proof}

We now turn to $k$\nobreakdash-Cayley regressions with $k\geq4$.
By Theorem~\ref{dihedral-2-cayley-regression},
it suffices to consider prime-power values of $n$.
We first treat powers of $2$.
The case $n=2$ was settled in Lemma~\ref{lem:D4-Cayley}.
For $n=2^m$ with $m\geq3$, we construct
a $4$-valent generalized Cayley graph that is not vertex-transitive.

\begin{lem}\label{2-power-not-4-cayley}
	Let $G=D_{2n}$, where $n=2^m$ and $m\geq3$.
	Then $G$ is not a $4$\nobreakdash-Cayley regression.
\end{lem}

\begin{proof}
Let $\alpha$ be conjugation by $b$, that is, $\alpha(g)=bgb^{-1}$ for $g\in G$.
Then $\alpha$ is an automorphism of $G$ with $\alpha^2=\mathrm{id}$, since $b^2=1$.
As $ba^kb^{-1}=a^{-k}$, we have
\[
\alpha(a^k)=a^{-k}\quad\text{and}\quad\alpha(a^kb)=a^{-k}b\qquad\text{for every integer }k.
\]
We compute $\omega_\alpha(G)$.
For $g=a^k$, we have $\alpha(g^{-1})g=\alpha(a^{-k})a^k=a^{2k}$.
For $g=a^kb$, we have $g^{-1}=g$, and so
$\alpha(g^{-1})g=a^{-k}b\cdot a^kb=a^{-k}a^{-k}=a^{-2k}$.
Hence $\omega_\alpha(G)=\{a^{2k}\mid k\in\mathbb Z_n\}=\langle a^2\rangle$.

Let $S=\{a,b,ab,a^{-1}b\}$.
The element $a$ lies in $A$, while $b$, $ab$ and $a^{-1}b$ lie in $Ab$ and are pairwise
distinct because $1$, $a$ and $a^{-1}$ are pairwise distinct; thus $|S|=4$.
Since $n$ is even, $a\notin\langle a^2\rangle$, and the other three elements of $S$ lie
outside $A$; hence $S\cap\omega_\alpha(G)=\varnothing$.
Moreover,
\[
\alpha(S)=\{a^{-1},b,a^{-1}b,ab\}=S^{-1}.
\]
Thus $S$ is a generalized Cayley subset of size $4$.
Put $\Gamma=GC(G,S,\alpha)$.

For $k\in\mathbb Z_n$, put $x_k=a^k$ and $z_k=a^{-k}b$, so that
$G=\{x_k\mid k\in\mathbb Z_n\}\mathbin{\dot\cup}\{z_k\mid k\in\mathbb Z_n\}$.
Recall that the neighbors of a vertex $g$ are the elements of $\alpha(g)S$.
Since $\alpha(x_k)=a^{-k}$, the neighbors of $x_k$ are
\[
a^{-k}a=x_{1-k},\quad a^{-k}b=z_k,\quad a^{-k}ab=z_{k-1},\quad a^{-k}a^{-1}b=z_{k+1}.
\]
Since $\alpha(z_k)=a^kb$ and $ba^j=a^{-j}b$ for every $j$, the neighbors of $z_k$ are
\[
a^kb\cdot a=z_{1-k},\quad a^kb\cdot b=x_k,\quad a^kb\cdot ab=x_{k-1},\quad a^kb\cdot a^{-1}b=x_{k+1}.
\]
Consequently, the edges of $\Gamma$ are exactly
\[
\{x_k,x_{1-k}\},\quad \{z_k,z_{1-k}\},\quad\text{and}\quad
\{x_k,z_\ell\}\ \text{with}\ \ell-k\in\{0,\pm1\},
\]
where $k,\ell\in\mathbb Z_n$.
Since $n\geq8$, two integers $c$ and $d$ with $0<|c-d|\leq7$ are not congruent modulo $n$,
because $n$ does not divide $c-d$; we use this fact twice below.

We now compare the vertices $x_0$ and $x_2$.
The vertex $x_0$ lies in the triangle $\{x_0,x_1,z_0\}$: indeed,
$\{x_0,x_1\}$ is an edge of the first type with $k=0$, and $\{x_0,z_0\}$ and $\{x_1,z_0\}$
are edges of the third type with $\ell-k=0$ and $\ell-k=-1$, respectively.
We claim that $x_2$ lies in no triangle, that is, no two neighbors of $x_2$ are adjacent.
The neighbors of $x_2$ are $x_{-1}$, $z_1$, $z_2$ and $z_3$.
First, the neighbors of $x_{-1}$ of the form $z_\ell$ are $z_{-2}$, $z_{-1}$ and $z_0$.
None of them lies in $\{z_1,z_2,z_3\}$: for $c\in\{-2,-1,0\}$ and $d\in\{1,2,3\}$ we have
$1\leq d-c\leq5$, so $c\not\equiv d\pmod n$ by the fact above.
Second, two vertices $z_\ell$ and $z_{\ell'}$ are adjacent only if $\ell+\ell'\equiv1\pmod n$;
for distinct $\ell,\ell'\in\{1,2,3\}$, we have $2\leq(\ell+\ell')-1\leq4$, so
$\ell+\ell'\not\equiv1\pmod n$ by the fact above.
This proves the claim.

Since automorphisms of $\Gamma$ preserve triangles, no automorphism maps $x_0$ to $x_2$.
Thus $\Gamma$ is not vertex-transitive, and so $G$ is not a $4$\nobreakdash-Cayley regression.
\end{proof}

The remaining power of $2$ is $n=4$.
Before treating $D_8$, we record a family of automorphisms of dihedral groups
for which every generalized Cayley graph is a Cayley graph;
it will also be used for $D_{2n}$ with $n$ an odd prime power.

\begin{lem}\label{lem:omega-A}
	Let $G=D_{2n}$, where $n\geq2$, and
	let $\alpha\in\operatorname{Aut}(G)$ be defined by $\alpha(a)=a^{-1}$ and
	$\alpha(b)=a^jb$, where $n$ or $j$ is odd.
	Then $\omega_\alpha(G)=A$, and
	$GC(G,S,\alpha)\cong\operatorname{Cay}(G,S)$ for every generalized
	Cayley subset $S$ of $G$ induced by $\alpha$.
\end{lem}

\begin{proof}
Since $\alpha^2(a)=a$ and $\alpha^2(b)=\alpha(a^j)\alpha(b)=a^{-j}a^jb=b$, we have
$\alpha^2=\mathrm{id}$.
Moreover, $\alpha(a^t)=a^{-t}$ and $\alpha(a^tb)=a^{j-t}b$ for every integer $t$;
in particular, $\alpha(A)=A$ and $\alpha(Ab)=Ab$.

We first compute $\omega_\alpha(G)$.
For $g=a^t$, we have $\alpha(g^{-1})g=a^ta^t=a^{2t}$.
For $g=a^tb$, we have $g^{-1}=g$, and so
$\alpha(g^{-1})g=a^{j-t}b\cdot a^tb=a^{j-t}a^{-t}=a^{j-2t}$.
Thus $\omega_\alpha(G)=\{a^{2t}\mid t\in\mathbb Z_n\}\cup\{a^{j-2t}\mid t\in\mathbb Z_n\}\subseteq A$.
If $n$ is odd, then $2$ is invertible modulo $n$, so the elements $a^{2t}$ already
exhaust $A$.
If $n$ is even and $j$ is odd, then the elements $a^{2t}$ are exactly the even powers
of $a$ and the elements $a^{j-2t}$ are exactly the odd powers of $a$, so together they
exhaust $A$.
In either case, $\omega_\alpha(G)=A$.

Now let $S$ be a generalized Cayley subset of $G$ induced by $\alpha$, and put
$\Gamma=GC(G,S,\alpha)$.
Since $S\cap\omega_\alpha(G)=\varnothing$, we have $S\subseteq Ab$.
Since every element of $Ab$ is an involution, we have $S^{-1}=S$, and hence $\alpha(S)=S$.
As also $1\notin S$, the Cayley graph $\operatorname{Cay}(G,S)$ is defined.

We next observe that both graphs are bipartite with parts $A$ and $Ab$.
In $\Gamma$, the neighbors of $g$ are the elements of $\alpha(g)S$; since $\alpha$
preserves $A$ and $Ab$ and $S\subseteq Ab$, the set $\alpha(g)S$ lies in the coset of $A$
not containing $g$.
Similarly, in $\operatorname{Cay}(G,S)$ the neighbors of $g$ are the elements of $gS$,
which also lie in the other coset.
Hence it suffices to compare adjacency between vertices $g\in A$ and $h\in Ab$.

Define $\phi:G\to G$ by $\phi(g)=g^{-1}$.
This is a bijection which maps $A$ onto $A$ and $Ab$ onto $Ab$.
Let $g\in A$ and $h\in Ab$.
Then $\alpha(g^{-1})=g$ and $h^{-1}=h$, so
\[
\phi(g)^{-1}\phi(h)=gh=\alpha(g^{-1})h.
\]
Thus $\phi(g)$ and $\phi(h)$ are adjacent in $\operatorname{Cay}(G,S)$ if and only if
$g$ and $h$ are adjacent in $\Gamma$.
Therefore $\phi$ is an isomorphism from $\Gamma$ to $\operatorname{Cay}(G,S)$.

Since $S$ was arbitrary, this completes the proof.
\end{proof}

\begin{lem}\label{D8-Cayley}
	The dihedral group $D_8$ is a Cayley regression.
\end{lem}

\begin{proof}
Let $G=D_8$, so that $n=4$.
Every automorphism of $G$ maps $a$ to $a^{\pm1}$ and $b$ to some $a^jb$.
If $a\mapsto a$ and $b\mapsto a^jb$, then the automorphism is involutory if and only if
$a^{2j}b=b$, that is, $j\in\{0,2\}$, and $j=0$ gives the identity.
If $a\mapsto a^{-1}$ and $b\mapsto a^jb$, then $b\mapsto a^jb\mapsto a^{-j}a^jb=b$,
so the automorphism is involutory for every $j$.
Hence the involutory automorphisms of $G$ are
\[
\delta:\ a\mapsto a,\ b\mapsto a^2b
\qquad\text{and}\qquad
\alpha_j:\ a\mapsto a^{-1},\ b\mapsto a^jb\quad(0\leq j\leq3).
\]
Let $\beta\in\operatorname{Aut}(G)$ be given by $\beta(a)=a$ and $\beta(b)=ab$.
Evaluating on $a$ and $b$, one checks that $\beta\alpha_0\beta^{-1}=\alpha_2$ and
$\beta\alpha_1\beta^{-1}=\alpha_3$.
By Proposition~\ref{prop:conjugate-GC}, it therefore suffices to consider
$\alpha\in\{\delta,\alpha_0,\alpha_1\}$.
For $\alpha=\alpha_1$, every generalized Cayley graph induced by $\alpha$ is isomorphic
to a Cayley graph of $G$ by Lemma~\ref{lem:omega-A}, since $j=1$ is odd.

From now on, let $\alpha\in\{\delta,\alpha_0\}$, let $S$ be a generalized Cayley subset
induced by $\alpha$, and put $\Gamma=GC(G,S,\alpha)$.
For $\alpha=\delta$, we have $\delta(a^{-k})a^k=1$ and
$\delta(a^kb)a^kb=a^ka^2b\cdot a^kb=a^2$; for $\alpha=\alpha_0$, the computation in the
proof of Lemma~\ref{2-power-not-4-cayley} gives $\omega_{\alpha_0}(G)=\langle a^2\rangle$.
Hence
\[
\omega_\delta(G)=\omega_{\alpha_0}(G)=\{1,a^2\}.
\]
Put $U=G\setminus\{1,a^2\}$.
Then $S\subseteq U$, so $|S|\leq6$.
By Lemma~\ref{lem:dihedral-prime-power-3}, $D_8$ is a $3$\nobreakdash-Cayley regression,
so we may assume that $4\leq|S|\leq6$.

We first record a useful observation.
Let $H$ be a subgroup of index two in $G$; the subgroups of index two in $G$ are
\[
A,\quad \langle a^2,b\rangle=\{1,a^2,b,a^2b\},\quad \langle a^2,ab\rangle=\{1,a^2,ab,a^3b\}.
\]
Since $H$ has index two, it is normal in $G$, and $G=H\mathbin{\dot\cup}xH$ for any $x\notin H$;
each of the two cosets of $H$ has four elements.
The automorphisms $\delta$ and $\alpha_0$ are inner, namely conjugation by $a$ and by $b$,
respectively, so they map every normal subgroup onto itself; thus $\alpha(H)=H$.
Consequently, $\alpha$ induces an automorphism of the group $G/H$ of order two, which must be
the identity.
In other words, $\alpha(x)\in xH$ for every $x\in G$.

Now let $g,h\in G$.
Applying the previous paragraph to $x=g^{-1}$, we obtain $\alpha(g^{-1})\in g^{-1}H$, and since
$H$ is normal,
\[
\alpha(g^{-1})h\in g^{-1}Hh=g^{-1}hH.
\]
Thus $\alpha(g^{-1})h$ and $g^{-1}h$ lie in the same coset of $H$, and so
\[
\alpha(g^{-1})h\notin H\iff g^{-1}h\notin H.
\]
Note that $G\setminus H$ is a generalized Cayley subset induced by $\alpha$: indeed,
$\omega_\alpha(G)=\{1,a^2\}\subseteq H$, and $\alpha(G\setminus H)=G\setminus H=(G\setminus H)^{-1}$.
In $GC(G,G\setminus H,\alpha)$, two vertices $g$ and $h$ are adjacent if and only if
$\alpha(g^{-1})h\notin H$, while in $\operatorname{Cay}(G,G\setminus H)$ they are adjacent
if and only if $g^{-1}h\notin H$.
By the equivalence above, these two graphs have the same edges, and hence they coincide.
Finally, $g^{-1}h\notin H$ if and only if $gH\neq hH$, that is, $g$ and $h$ lie in different
cosets of $H$.
Since each coset has four elements, every vertex is adjacent to exactly the four vertices of the
other coset and to no vertex of its own coset.
Therefore
\begin{equation}\tag{$\ast$}\label{eq:index-two-complement}
GC(G,G\setminus H,\alpha)=\operatorname{Cay}(G,G\setminus H)\cong K_{4,4}.
\end{equation}

\emph{Case 1: $|S|=4$.}
Since $\alpha(S)=S^{-1}$, the set $S$ is invariant under the involution
$x\mapsto\alpha(x^{-1})$ of $U$.

Suppose first that $\alpha=\delta$.
The involution $x\mapsto\delta(x^{-1})$ has no fixed points on $U$, and its orbits are
$\{a,a^3\}$, $\{b,a^2b\}$ and $\{ab,a^3b\}$.
Hence $S$ is the union of two of these pairs.
Adding $\{1,a^2\}$ to the remaining pair gives one of the three subgroups $H$ of index
two listed above, and $S=G\setminus H$.
Thus $\Gamma\cong K_{4,4}$ by~\eqref{eq:index-two-complement}.

Suppose next that $\alpha=\alpha_0$.
The involution $x\mapsto\alpha_0(x^{-1})$ fixes $a$, $a^3$, $b$ and $a^2b$, and
interchanges $ab$ and $a^3b$.
Hence either $S=\{a,a^3,b,a^2b\}$, or $S=\{ab,a^3b,x,y\}$ for distinct
$x,y\in\{a,a^3,b,a^2b\}$; this gives seven possibilities.
If $S$ is $\{a,a^3,b,a^2b\}$, $\{a,a^3,ab,a^3b\}$ or $Ab$, then $S$ is the complement of
$\langle a^2,ab\rangle$, $\langle a^2,b\rangle$ or $A$, respectively,
and $\Gamma\cong K_{4,4}$ by~\eqref{eq:index-two-complement}.
In the remaining four cases, $\{x,y\}$ consists of one element of $\{a,a^3\}$
and one element of $\{b,a^2b\}$.
Let $\tau\in\operatorname{Aut}(G)$ be given by $\tau(a)=a$ and $\tau(b)=a^2b$, and put
$S_0=\{a,b,ab,a^3b\}$.
The automorphisms $\alpha_0$, $\tau$ and $\alpha_0\tau$ commute with $\alpha_0$, and
\[
\begin{aligned}
\alpha_0(S_0)&=\{a^3,b,ab,a^3b\},\\
\tau(S_0)&=\{a,a^2b,ab,a^3b\},\\
\alpha_0\tau(S_0)&=\{a^3,a^2b,ab,a^3b\}
\end{aligned}
\]
are the other three such subsets.
For $\gamma\in\{\alpha_0,\tau,\alpha_0\tau\}$, Proposition~\ref{prop:conjugate-GC} gives
\[
GC(G,S_0,\alpha_0)\cong GC(G,\gamma(S_0),\gamma\alpha_0\gamma^{-1})=GC(G,\gamma(S_0),\alpha_0),
\]
so it suffices to consider $S=S_0$.
Writing $N(g)=\alpha_0(g)S_0$ for the neighborhood of a vertex $g$, we obtain
\[
\begin{aligned}
N(1)&=\{a,b,ab,a^3b\}, & N(a^2)&=\{a^3,ab,a^2b,a^3b\},\\
N(a)&=\{1,b,a^2b,a^3b\}, & N(a^3)&=\{a^2,b,ab,a^2b\},\\
N(b)&=\{1,a,a^3,a^3b\}, & N(ab)&=\{1,a^2,a^3,a^2b\},\\
N(a^3b)&=\{1,a,a^2,b\}, & N(a^2b)&=\{a,a^2,a^3,ab\}.
\end{aligned}
\]
Hence $X_1=\{1,a,b,a^3b\}$ and $X_2=\{a^2,a^3,ab,a^2b\}$ both induce $K_4$,
and the only edges between them are
$\{1,ab\}$, $\{a,a^2b\}$, $\{b,a^3\}$ and $\{a^3b,a^2\}$,
which form a perfect matching.
Therefore $\Gamma\cong K_4\square K_2$.
This is isomorphic to $\operatorname{Cay}(G,\{a,a^2,a^3,b\})$, in which $A$ and $Ab$ induce
copies of $K_4$ and the edges $\{g,gb\}$ form a perfect matching between them.

\emph{Case 2: $|S|\in\{5,6\}$.}
Since $\{1,a^2\}$ is invariant under both $\alpha$ and inversion, we have
$\alpha(U)=U=U^{-1}$ and $U\cap\omega_\alpha(G)=\varnothing$.
Thus $U$ is a generalized Cayley subset, and $\Delta=GC(G,U,\alpha)$ is a $6$-valent
graph on eight vertices.
Hence $\overline{\Delta}$ is $1$-valent, that is, a perfect matching.
If $|S|=6$, then $S=U$ and $\Gamma=\Delta$; thus $\Gamma$ is the complement of a perfect
matching on eight vertices, and so
$\Gamma\cong\overline{4K_2}\cong\operatorname{Cay}(G,G\setminus\{1,a^2\})$.

If $|S|=5$, write $U\setminus S=\{t\}$.
From $\alpha(U)=U^{-1}$ and $\alpha(S)=S^{-1}$, we obtain $\alpha(t)=t^{-1}$.
Hence $\{t\}$ is a generalized Cayley subset, and $GC(G,\{t\},\alpha)$ is $1$-valent,
that is, a perfect matching.
Since $U=S\mathbin{\dot\cup}\{t\}$, the edge set of $\Delta$ is the disjoint union of the
edge sets of $\Gamma$ and $GC(G,\{t\},\alpha)$.
Consequently, $\overline{\Gamma}$ is the union of the two edge-disjoint perfect matchings
$\overline{\Delta}$ and $GC(G,\{t\},\alpha)$.
Such a union is a $2$-valent graph whose cycles alternate between the two matchings,
so all its cycles have even length.
On eight vertices, this leaves
$\overline{\Gamma}\cong C_8\cong\operatorname{Cay}(G,\{b,ab\})$ or
$\overline{\Gamma}\cong 2C_4\cong\operatorname{Cay}(G,\{a,a^{-1}\})$.
Since the complement of a Cayley graph of $G$ is again a Cayley graph of $G$,
$\Gamma$ is isomorphic to a Cayley graph of $G$.

Hence, in all cases, $\Gamma$ is isomorphic to a Cayley graph of $G$, so $D_8$ is a
Cayley regression.
\end{proof}

Finally, we show that $D_{2n}$ is a Cayley regression whenever $n$ is an odd prime power.

\begin{lem}\label{odd-prime-power-cayley-regression}
	Let $G=D_{2p^m}$, where $p$ is an odd prime and $m\geq1$.
	Then $G$ is a Cayley regression.
\end{lem}

\begin{proof}
Put $n=p^m$; then $n\geq3$ is odd.
Let $\alpha$ be an involutory automorphism of $G$.
Since $n\geq3$ and every element of $Ab$ is an involution, every element of order $n$
lies in $A$; hence $A$ is the unique cyclic subgroup of order $n$ in $G$,
and so $\alpha(A)=A$ and $\alpha(Ab)=Ab$.
Thus there exist $u\in\mathbb Z_n^{*}$ and $v\in\mathbb Z_n$ such that
$\alpha(a)=a^u$ and $\alpha(b)=a^vb$.
Then
\[
\alpha^2(a)=a^{u^2}\quad\text{and}\quad\alpha^2(b)=\alpha(a^v)\alpha(b)=a^{uv}a^vb=a^{v(u+1)}b,
\]
so $\alpha^2=\mathrm{id}$ gives $u^2\equiv1\pmod n$ and $v(u+1)\equiv0\pmod n$.

The congruence $u^2\equiv1\pmod n$ means that $(u-1)(u+1)$ is divisible by $p^m$.
Since $u-1$ and $u+1$ differ by $2$ and $p$ is odd, $p$ does not divide both of them,
and hence $p^m$ divides one of them; that is, $u\equiv\pm1\pmod n$.
If $u\equiv1\pmod n$, then $2v\equiv0\pmod n$, and since $n$ is odd, $v\equiv0\pmod n$;
thus $\alpha=\mathrm{id}$, and every generalized Cayley graph induced by $\alpha$ is a
Cayley graph of $G$.
If $u\equiv-1\pmod n$, then $\alpha(a)=a^{-1}$ and $\alpha(b)=a^vb$, and since $n$ is odd,
Lemma~\ref{lem:omega-A} shows that every generalized Cayley graph induced by $\alpha$ is
isomorphic to a Cayley graph of $G$.

Since $\alpha$ was an arbitrary involutory automorphism, $G$ is a Cayley regression.
\end{proof}

We can now complete the classification for valency at least four.

\begin{proof}[Proof of Theorem~\ref{thm:dihedral-k-cayley-regression}]
	Let $G=D_{2n}$, where $n\geq2$, and let $k\geq4$.
	Suppose first that $G$ is a $k$\nobreakdash-Cayley regression.
	Then $G$ is a $2$\nobreakdash-Cayley regression, and hence $n$ is a prime power
	by Theorem~\ref{dihedral-2-cayley-regression}.
	Moreover, $G$ is a $4$\nobreakdash-Cayley regression, so
	Lemma~\ref{2-power-not-4-cayley} excludes $n=2^m$ with $m\geq3$.
	Therefore $n$ is an odd prime power or $n\in\{2,4\}$.

	Conversely, $G$ is a Cayley regression, and hence a
	$k$\nobreakdash-Cayley regression, by
	Lemma~\ref{odd-prime-power-cayley-regression} if $n$ is an odd prime power,
	by Lemma~\ref{lem:D4-Cayley} if $n=2$, and
	by Lemma~\ref{D8-Cayley} if $n=4$.
\end{proof}

\begin{proof}[Proof of Corollary~\ref{classification-Cayley-regression}]
	If $D_{2n}$ is a Cayley regression, then it is in particular a
	$4$\nobreakdash-Cayley regression, and so $n$ is an odd prime power or $n\in\{2,4\}$
	by Theorem~\ref{thm:dihedral-k-cayley-regression}.
	Conversely, suppose that $n$ is an odd prime power or $n\in\{2,4\}$, and put
	$k=\max\{4,2n-1\}$.
	By Theorem~\ref{thm:dihedral-k-cayley-regression}, $D_{2n}$ is a
	$k$\nobreakdash-Cayley regression.
	Since every generalized Cayley graph $GC(D_{2n},S,\alpha)$ has valency
	$|S|\leq2n-1\leq k$, it follows that $D_{2n}$ is a Cayley regression.
\end{proof}

\section{Cayley regressions for some metacyclic groups}\label{sec:index-two}

In this section, we prove Theorems~\ref{thm:index-two-classification}
and~\ref{thm:split-metacyclic-odd}.
In both cases, $G$ has a cyclic normal subgroup $A$ with $G/A$ cyclic,
and we construct a $2$-valent generalized Cayley graph of $G$ that is not
isomorphic to a Cayley graph of $G$.

\subsection{Groups with a cyclic subgroup of index two}\label{subsec:index-two}

The cyclic case is known from \cite{LiaoLiuTang}, and the dihedral case was settled in
Section~\ref{sec:dihedral}.
It remains to show that groups with a cyclic subgroup of index two that are neither
cyclic nor dihedral are not $2$\nobreakdash-Cayley regressions.

Let $G$ be a group with a cyclic subgroup $A=\langle a\rangle$ of index two,
where $|A|=m$, and let $b\in G\setminus A$.
Then $A\trianglelefteq G$, $G=A\mathbin{\dot\cup}Ab$, and there are integers $t$ and $r$
with $\gcd(r,m)=1$ such that $b^2=a^t$ and $b^{-1}ab=a^r$.
The relations $a^m=1$, $b^2=a^t$ and $b^{-1}ab=a^r$ reduce every word in $a,b$ to the
form $a^ib^j$ with $0\leq i<m$ and $j\in\{0,1\}$, so the group they define has order at most $2m$;
since it maps onto $G$ and $|G|=2m$, these relations form a presentation of $G$.
Since conjugation by $b^2\in A$ is trivial on $A$, we have $r^2\equiv1\pmod m$;
since $b$ commutes with $b^2=a^t$, we have $a^t=b^{-1}a^tb=a^{rt}$, that is,
$(r-1)t\equiv0\pmod m$.

\begin{proof}[Proof of Theorem~\ref{thm:index-two-classification}\ref{idx:i}]
Retain the notation above, and suppose that $G$ is neither cyclic nor dihedral.
Since $|G|=2m$ and every group of order two is cyclic, our assumption that $G$ is not
cyclic forces $m>1$; that is, $a\neq1$.

We first construct an involutory automorphism of $G$.
Since $b^2\in A$ commutes with $a$, we have
$bab^{-1}=b^{-1}(b^2ab^{-2})b=b^{-1}ab=a^r$.
Define $\alpha$ on the generators by $\alpha(a)=a^{-r}$ and $\alpha(b)=b^{-1}$.
These images satisfy the defining relations of $G$:
\begin{align*}
\alpha(a)^m&=a^{-rm}=1,\\
\alpha(b)^2&=b^{-2}=a^{-t}=a^{-rt}=\alpha(a)^t,\\
\alpha(b)^{-1}\alpha(a)\alpha(b)&=ba^{-r}b^{-1}=a^{-r^2}=\alpha(a)^r,
\end{align*}
where the second line uses $(r-1)t\equiv0\pmod m$.
Hence $\alpha$ extends to an endomorphism of $G$.
Since $\alpha^2(a)=\alpha(a)^{-r}=a^{r^2}=a$ and $\alpha^2(b)=\alpha(b)^{-1}=b$,
we have $\alpha^2=\mathrm{id}$; in particular, $\alpha$ is an involutory automorphism of $G$.
Moreover, $\alpha(A)=A$, so $\alpha$ induces the identity on $G/A$, which has order two;
that is, $\alpha(g)$ and $g$ lie in the same coset of $A$ for every $g\in G$.

Let $S=\{b,ab\}$.
Since $a\neq1$, we have $|S|=2$, and $S\subseteq Ab$.
By the previous paragraph, $\alpha(g^{-1})g\in A$ for every $g\in G$, so
$\omega_\alpha(G)\subseteq A$ and hence $S\cap\omega_\alpha(G)=\varnothing$.
Moreover, $\alpha(b)=b^{-1}$, and since $b^{-1}a^{-1}=a^{-r}b^{-1}$, we have
$\alpha(ab)=a^{-r}b^{-1}=b^{-1}a^{-1}=(ab)^{-1}$.
Thus $\alpha(S)=S^{-1}$, and $S$ is a generalized Cayley subset of size two.
Put $\Gamma=GC(G,S,\alpha)$; then $\Gamma$ is $2$-valent, and the neighbors of a vertex $x$
are $\alpha(x)b$ and $\alpha(x)ab$.

We show that $\Gamma$ is connected.
Let $x\in G$.
The vertex $x$ is adjacent to $y=\alpha(x)b$, and $y$ is adjacent to
\[
\alpha(y)ab=\alpha^2(x)\alpha(b)ab=xb^{-1}ab=xa^r.
\]
Hence $x$ and $xa^r$ lie in the same connected component of $\Gamma$.
Applying this repeatedly, all the vertices $x,xa^r,xa^{2r},\ldots$ lie in this component.
Since $\gcd(r,m)=1$, the element $a^r$ generates $A$, so the whole coset $xA$ lies in a
single component.
Finally, $y=\alpha(x)b$ lies in the coset of $A$ different from $xA$, because $\alpha(x)\in xA$.
Thus both cosets of $A$ lie in the same component, and $\Gamma$ is connected.

If $\Gamma$ were isomorphic to a Cayley graph of $G$, then $G$ would admit a connected
$2$-valent Cayley graph, and so $G$ would be cyclic or dihedral by
Lemma~\ref{lem:connected-2-valent-Cayley}, contrary to our assumption.
Hence $G$ is not a $2$\nobreakdash-Cayley regression.
\end{proof}

\begin{proof}[Proof of Theorem~\ref{thm:index-two-classification}\ref{idx:ii}]
Let $G$ be a group with a cyclic subgroup of index two.
We consider three cases: $G$ is neither cyclic nor dihedral, $G$ is cyclic, or $G$ is dihedral.

If $G$ is neither cyclic nor dihedral, then $G$ is not a $2$\nobreakdash-Cayley regression by
Theorem~\ref{thm:index-two-classification}\ref{idx:i}.
Since every Cayley regression is a $2$\nobreakdash-Cayley regression, $G$ is not a Cayley regression.

If $G$ is cyclic, then $|G|$ is even, because $G$ has a subgroup of index two.
By~\cite[Theorem~2.11]{LiaoLiuTang}, a cyclic group of order $n>1$ is a Cayley regression
if and only if $n$ is an odd prime power or $n\in\{2,4,6,10\}$.
As $|G|$ is even, it follows that $G$ is a Cayley regression if and only if
$G$ is isomorphic to $C_2$, $C_4$, $C_6$ or $C_{10}$.

If $G\cong D_{2n}$ with $n\geq2$, then, by Corollary~\ref{classification-Cayley-regression},
$G$ is a Cayley regression if and only if $n$ is an odd prime power or $n\in\{2,4\}$.

Combining the three cases proves part~\ref{idx:ii}.
\end{proof}

For $n\geq2$, consider the dicyclic group $T_{4n}$ of order $4n$ and the group $SD_{8n}$ of order $8n$ given by
\begin{align*}
T_{4n}&=\langle a,b\mid a^{2n}=b^4=1,\ a^n=b^2,\ b^{-1}ab=a^{-1}\rangle,\\
SD_{8n}&=\langle a,b\mid a^{4n}=b^2=1,\ ba=a^{2n-1}b\rangle.
\end{align*}
When $n$ is a power of $2$, $SD_{8n}$ is the semidihedral group of order $8n$.

\begin{cor}\label{cor:T4n}
	For $n\geq2$, neither $T_{4n}$ nor $SD_{8n}$ is a $2$\nobreakdash-Cayley regression.
\end{cor}

\begin{proof}
Let $G$ be $T_{4n}$ or $SD_{8n}$.
We check that $G$ satisfies the hypotheses of
Theorem~\ref{thm:index-two-classification}\ref{idx:i}.

\emph{$G$ has a cyclic subgroup of index two.}
The subgroup $\langle a\rangle$ is cyclic of order $N$, where $N=2n$ if $G=T_{4n}$ and
$N=4n$ if $G=SD_{8n}$.
Since $|G|=2N$, it has index two.
Note that $N\geq4$.

\emph{$G$ is not cyclic.}
A cyclic group is abelian, but $G$ is not:
in $T_{4n}$ we have $b^{-1}ab=a^{-1}\neq a$, because $a$ has order $2n\geq4$;
in $SD_{8n}$ we have $ba=a^{2n-1}b\neq ab$, because $a$ has order $4n$.

\emph{$G$ is not dihedral.}
Recall that in the dihedral group $D_{2N}$ with $N\geq3$, the rotations form the unique
cyclic subgroup of order $N$, and every element outside it (a reflection) is an involution.
So if $G$ were isomorphic to $D_{2N}$, then $\langle a\rangle$ would correspond to the
rotations, and every element of $G$ outside $\langle a\rangle$ would be an involution.
This is not the case.
In $T_{4n}$, the element $b$ lies outside $\langle a\rangle$, but $b^2=a^n\neq1$.
In $SD_{8n}$, the element $ab$ lies outside $\langle a\rangle$, and using $b^2=1$ and
$bab^{-1}=a^{2n-1}$ we get
\[
(ab)^2=a\,(bab^{-1})\,b^2=a\cdot a^{2n-1}=a^{2n}\neq1.
\]
Hence $G$ is not dihedral.

Therefore $G$ is not a $2$\nobreakdash-Cayley regression by
Theorem~\ref{thm:index-two-classification}\ref{idx:i}.
\end{proof}

\subsection{Split metacyclic groups with a cyclic normal subgroup of odd index}\label{subsec:odd-index}

We now prove the analogous result for split metacyclic groups $C_m\rtimes C_q$ with $q$ odd.

\begin{proof}[Proof of Theorem~\ref{thm:split-metacyclic-odd}]
Let $G=C_m\rtimes C_q$, where $m>2$ and $q\geq3$ is odd.
Let $A=\langle a\rangle\cong C_m$ be the normal factor and $B=\langle b\rangle\cong C_q$
a complement, so that $G=AB$ and $A\cap B=\{1\}$.
Conjugation by $b$ induces an automorphism of $A$, so $b^{-1}ab=a^r$ for some integer $r$
with $\gcd(r,m)=1$, and $r^q\equiv1\pmod m$ because $b^q=1$.
Thus
\[
G=\langle a,b\mid a^m=b^q=1,\ b^{-1}ab=a^r\rangle,
\]
and every element of $G$ has a unique expression $a^ib^j$ with $i\in\mathbb Z_m$ and
$j\in\mathbb Z_q$.

We first construct an involutory automorphism of $G$.
Define $\alpha$ on the generators by $\alpha(a)=a^{-1}$ and $\alpha(b)=b$.
Then $\alpha(a)^m=\alpha(b)^q=1$ and
\[
\alpha(b)^{-1}\alpha(a)\alpha(b)=b^{-1}a^{-1}b=(b^{-1}ab)^{-1}=a^{-r}=\alpha(a)^r,
\]
so the images satisfy the defining relations, and $\alpha$ extends to an endomorphism of $G$.
Since $\alpha^2(a)=a$ and $\alpha^2(b)=b$, we have $\alpha^2=\mathrm{id}$, and $\alpha$ is
an involutory automorphism of $G$; explicitly, $\alpha(a^ib^j)=a^{-i}b^j$.

Let $S=\{b,b^{-1}\}$.
Since $b$ has order $q\geq3$, we have $b\neq b^{-1}$, so $|S|=2$.
Since $\alpha(A)=A$ and $\alpha(b)=b$, the automorphism $\alpha$ induces the identity on
$G/A$; hence $\alpha(g^{-1})\in g^{-1}A$, and so $\alpha(g^{-1})g\in g^{-1}Ag=A$ for every
$g\in G$.
Thus $\omega_\alpha(G)\subseteq A$.
As $b^{\pm1}\notin A$, we get $S\cap\omega_\alpha(G)=\varnothing$.
Moreover, $\alpha(S)=S=S^{-1}$.
Therefore $S$ is a generalized Cayley subset of size two, and we put $\Gamma=GC(G,S,\alpha)$.

Since $\alpha(a^ib^j)=a^{-i}b^j$, the neighbors of $a^ib^j$ are
$\alpha(a^ib^j)b=a^{-i}b^{j+1}$ and $\alpha(a^ib^j)b^{-1}=a^{-i}b^{j-1}$.
Let $\varphi$ be the permutation of $G$ given by $\varphi(a^ib^j)=a^{-i}b^{j+1}$.
Its inverse is $a^ib^j\mapsto a^{-i}b^{j-1}$, so the two neighbors of each vertex $x$
are $\varphi(x)$ and $\varphi^{-1}(x)$.
Hence the connected components of $\Gamma$ are precisely the orbits of
$\langle\varphi\rangle$ on $G$, and the orbit of $x$ has length equal to the least
positive integer $k$ with $\varphi^k(x)=x$.

By induction on $k$, we have $\varphi^k(a^ib^j)=a^{(-1)^ki}b^{j+k}$.
By the uniqueness of the expression $a^ib^j$, it follows that $\varphi^k(a^ib^j)=a^ib^j$
if and only if $j+k\equiv j\pmod q$ and $(-1)^ki\equiv i\pmod m$, that is,
if and only if $k$ is divisible by $q$ and $(-1)^ki\equiv i\pmod m$.
For the vertex $1=a^0b^0$, the second condition holds for every $k$,
so the orbit of $1$ has length $q$.
For the vertex $a=a^1b^0$, since $m>2$ we have $-1\not\equiv1\pmod m$, so the second
condition holds if and only if $k$ is even.
Thus $\varphi^k(a)=a$ if and only if $k$ is divisible by both $2$ and $q$; as $q$ is odd,
this means that $k$ is divisible by $2q$, and so the orbit of $a$ has length $2q$.

Therefore $\Gamma$ has connected components of orders $q$ and $2q$.
Since all connected components of a vertex-transitive graph are isomorphic,
$\Gamma$ is not vertex-transitive, and so $G$ is not a $2$\nobreakdash-Cayley regression.
\end{proof}

 \section*{Funding information}
This research was supported by the National Natural Science Foundation of China
(No.~12561002) and by the Hainan Provincial Natural Science Foundation of China
(No.~126MS0135).

\end{document}